\documentclass[11pt]{amsart}
\usepackage[margin=1.12in]{geometry}
\usepackage{amsmath,amssymb,amsthm,mathtools,booktabs,array,longtable}
\usepackage{microtype}
\usepackage{xurl}
\usepackage[hidelinks]{hyperref}
\hypersetup{pdftitle={Rigidity of complete self-shrinkers}}
\numberwithin{equation}{section}
\newtheorem{theorem}{Theorem}[section]
\newtheorem{proposition}{Proposition}[section]

\newtheorem{lemma}[theorem]{Lemma}
\newtheorem{example}{Example}[section]
\newtheorem{corollary}[theorem]{Corollary}
\theoremstyle{remark}\newtheorem{remark}[theorem]{Remark}
\newcommand{\R}{\mathbb R}
\newcommand{\Lop}{\mathcal L}
\newcommand{\tr}{\operatorname{tr}}
\allowdisplaybreaks[2]
\title[Rigidity of complete self-shrinkers]{Rigidity of complete self-shrinkers}
\author [Q. -M. Cheng, F. Li  and G. Wei]{Qing-Ming Cheng, Fengjiang Li{$^+$}  and Guoxin Wei{$^+$}}
\address{Qing-Ming Cheng \newline
\indent
Mathematical Science Research Center, \newline 
\indent Chongqing University of Technology, Chongqing 400054, P. R. China. 
\vskip 1mm
\indent
Osaka Central  Advanced Mathematical Institute \newline
\indent
Osaka Metropolitan University, Osaka 558-8585, Japan. \newline
\indent qingmingcheng@yahoo.com, chengqingming@cqut.edu.cn }
\address{Fengjiang Li \newline
\indent
Mathematical Science Research Center, \newline 
\indent Chongqing University of Technology, Chongqing 400054, P. R. China. \newline 
\indent fengjiangli@cqut.edu.cn }

\address{Guoxin Wei \newline
\indent   School of Mathematical Sciences, \newline 
\indent  South China Normal University, Guangzhou 510631, P. R. China. \newline 
\indent weiguoxin@tsinghua.org.cn}

\thanks{$+$ Fengjiang Li  and Guoxin Wei are corresponding authors}
\thanks{*This work was supported by  National Natural Science Foundation of China (Grant Nos.
12301062, 12571050, 12671065), Natural Science Foundation of Chongqing
(No.CSTB2024NSCQ-MSX0537), Japan Society for the Promotion of Science Grant-in-Aid
for Scientific Research (C) (Grant No. 25K06992).}
\date{November 2026}

\DeclareMathOperator{\Div}{div}
\theoremstyle{definition}\newtheorem{definition}[theorem]{Definition}

\begin{document}
\begin{abstract}
Let $X:M^n\to\R^{n+1}$ be an $n$-dimensional  complete self-shrinker
with constant squared norm $S$ of the second fundamental form. 
For $n=4$,  we prove that  if $f_3=\sum_i\lambda_i^3$ is constant, 
 where $\lambda_i$ denotes principal curvatures, 
then $S=0$, or $S=1$
and the image $X(M^n)$ in the rigid cases is one of  the hyperplane, the sphere  
or  the standard generalized cylinders. 
In every dimension, under assumption of sub-exponential growth 
we  prove the same  conclusion. 
Without assumption of sub-exponential growth, we prove  
if $S\leq 3+1/\sqrt n$, then the same  conclusion holds.

\end{abstract}
\maketitle

\section{Introduction and main results}

In this paper, we consider connected, complete hypersurfaces without boundary in Euclidean space. 
A hypersurface $X:M^n\to\R^{n+1}$ is called a self-shrinker if 
\begin{equation}\label{shrinker}
 H+\langle X,N\rangle=0,
\end{equation}
where  $N$ is a unit normal and $H$ is the mean curvature. \\
It is well-known that self-shrinkers arise as self-similar solutions and as models for singularities
of the mean curvature flow. Huisken's work \cite{Huisken} and the subsequent theory of Colding and Minicozzi \cite{CM} 
place self-shrinkers at the center of the study of mean curvature flow singularities.
The standard examples
are the hyperplane through the origin, the sphere $S^n(\sqrt n)$, and
the generalized cylinders $S^k(\sqrt k)\times\R^{n-k}$.

Let $\vec h$ be the second fundamental form, and let
$\lambda_1,\ldots,\lambda_n$ be the principal curvatures. We write
\begin{equation}\label{intro:traces}
 S=|\vec h|^2=\sum_i\lambda_i^2,\qquad
 f_r=\sum_i\lambda_i^r.
\end{equation}
The classification of self-shrinkers with constant $S$ is closely related
to gap phenomena for these curvature invariants. Cao--Li \cite{CaoLi}
established a first gap theorem under assumption of polynomial volume growth.
Cheng--Peng \cite{CP} developed a generalized maximum principle for the
$\Lop$-operator and obtained rigidity results without this volume growth
assumption. For second gap and pinching results, we refer to
Ding--Xin \cite{DX}, Cheng--Wei \cite{CW} and
Lei--Xu--Xu \cite{LeiXuXu}. In particular, under constant $S$ and polynomial
volume growth, Cheng--Wei \cite{CW} excluded the interval
$1<S<10/7$. 
Recently, Cheng--Li--Wei \cite{CLW3} have proven that 
for an $n$-dimensional  complete self-shrinker $X:M^n\to\R^{n+1}$
with constant squared norm of second fundamental form, up to an orthogonal
transformation, its image $X(M)$ is one of
\[
\R^n,\qquad S^n(\sqrt n),\qquad
S^k(\sqrt k)\times\R^{n-k}\quad(1\le k\le n-1),
\]
if scalar curvature $R$ satisfies $R\geq -\frac32$.
\\
Under the additional assumption that $f_3$ is constant,
Cheng--Wei--Yano \cite{CWY} obtained a gap up to $S<1.83379$ without
assuming polynomial volume growth. 
\\
In low dimensions, Cheng--Ogata \cite{CO} classified complete
$2$-dimensional self-shrinkers with constant $S$.
Cheng--Li--Wei classified complete $3$-dimensional self-shrinkers under
constant $S,f_4$ in \cite{CLWf4}, and under constant $S,\ f_3$ in
\cite{CLWf3}. Recent preprints of Fu--Yang \cite{FY} and Guan--Zhao \cite{GZ}
address the classification in dimension three under constant $S$ alone.\\
Recently, Cheng--Li \cite{cl26} has given a complete classification of complete 
self-shrinkers with constant $S, \ f_3$ and polynomial
volume growth. It is known that Ding--Xin \cite{DXVolume} and Cheng--Zhou \cite{ChengZhou}
has proven that for  $n$-dimensional complete self-shrinkers, polynomial 
volume growth and properness is equivalent.  On the other hand, from the following example, we know that 
 the condition of polynomial
volume growth is a strong condition.

\begin{example}\label{rem:cover}
The  complete immersed self-shrinker
\[
 X:\R\times\R^{n-1}\longrightarrow\R^{n+1},
 \qquad X(t,z)=(\cos t,\sin t,z)
\]
has $S=1$ and constant $f_3$, and its induced metric is Euclidean.
Its volumes of intrinsic balls are $\omega_n r^n$, while the inverse image
of $S^1(1)\times\{0\}$ is noncompact. Thus it is a complete nonproper
immersion.
This example also distinguishes intrinsic volume growth from extrinsic volume growth
counted with multiplicity.
\end{example}

The present paper concerns complete
self-shrinkers with constant $S$ and $f_3$. We use the curvature identities
and notation of Cheng--Li--Wei \cite{CLW1}.
Our first result gives the following:

\begin{theorem}\label{four-dimensional-classification}
Let $X:M^4\to\R^5$ be a 4-dimensional complete self-shrinker.
If $S$ and $f_3$ are constant, then $S=0$ or $S=1$.
Up to an orthogonal transformation, its image is a hyperplane through
the origin, $S^4(2)$, or $S^k(\sqrt k)\times\R^{4-k}$ for
$1\le k\le3$.
\end{theorem}
\begin{remark}
In our theorem, we remove the assumption of polynomial volume growth
in \cite{cl26}. 
\end{remark}

Since polynomial volume growth provides a different approach for studying
complete self-shrinkers  in arbitrary dimension as Colding and Minicozzi \cite{CM}, Cao and Li \cite{CaoLi} 
and so on, in the next theorem, we  weak polynomial volume growth to
sub-exponential volume growth.
\\
For a fixed $p\in M$, the intrinsic volume $ V_p(r)$ of a geodesic ball $B_r^M(p)$ of radius $r$ is defined by
\begin{equation}\label{eq:ball-volumes}
 V_p(r)=\operatorname{Vol}(B_r^M(p)),\qquad
 W_p(r)=\int_{B_r^M(p)}e^{-|X|^2/2}\,d\mu.
\end{equation}
Sub-exponential volume growth means 
$$
 \limsup_{r\to\infty}\frac{\log V_p(r)}r=0.
 $$

\begin{theorem}\label{thm:maingrowth}
Let $X:M^n\to\R^{n+1}$ be an $n$-dimensional  complete self-shrinker.
Suppose that $S$ and $f_3$ are constant and that $M$ has sub-exponential
volume growth, then $S=0$ or $S=1$, and, up to an orthogonal
transformation, the image $X(M)$ is one of
\[
\R^n,\qquad S^n(\sqrt n),\qquad
S^k(\sqrt k)\times\R^{n-k}\quad(1\le k\le n-1),
\]
where $\mathbb R^n$ denotes a hyperplane through the origin.
The same conclusion holds under the weaker assumption
\begin{equation}\label{eq:main-weighted-growth}
 \limsup_{r\to\infty}\frac{\log W_p(r)}r=0.
\end{equation}
\end{theorem}
\begin{remark}
For complete self-shrinkers, polynomial volume growth yields sub-exponential volume growth.
\end{remark}

For general dimension and without volume condition, we have the following second gap:
\begin{theorem}\label{thm:f3}
Let $X:M^n\to\R^{n+1}$ be an $n$-dimensional  complete self-shrinker.
Suppose that $S$ and $f_3$ are constant. If $S>1$, then
\begin{equation}\label{eq:main-second-gap}
 S>3+\frac1{\sqrt n}.
\end{equation}
Consequently, if $S\le3+1/\sqrt n$, then $S=0$ or $S=1$.
In this case, up to an orthogonal transformation, the image $X(M)$ is one of
\[
\R^n,\qquad S^n(\sqrt n),\qquad
S^k(\sqrt k)\times\R^{n-k}\quad(1\le k\le n-1).
\]
\end{theorem}

In Section~\ref{prelim}, we recall the known curvature estimates,
maximum principles and volume-growth conditions.
In Section~\ref{sec:algebra-cubic}, we compare the maximum sequence estimate
for $f_4$ with the minimum sequence identities for $H$ to prove
Theorem~\ref{thm:f3}. In Section~\ref{sec:four-dimensional}, we prove Theorem 1.1  by analyzing limiting principal
curvatures at infinity. Finally, in Section 5, we prove
Theorem~\ref{thm:maingrowth} using a weighted divergence identity.

\section{Preliminaries and notation}\label{prelim}

In this section, we recall the notation and basic formulas of
Cheng--Li--Wei \cite{CLW1}. We also recall maximum principles
and give the volume growth conditions.

\subsection{Notation, basic identities and estimates}\label{subsec:basic}
Let $X:M^n\to\R^{n+1}$ be an $n$-dimensional  self-shrinker, and let
$N=e_{n+1}$ be a local unit normal. Complete it to an orthonormal
frame $\{e_1,\ldots,e_n,e_{n+1}\}$ along the immersion, with dual
forms $\omega_i$, $1\le i\le n$. All Latin indices below range
from $1$ to $n$. Unless otherwise indicated, sums over indices are
ordered. Our conventions are
\[
dX=\sum_i\omega_i e_i,\qquad
de_i=\sum_j\omega_{ij}e_j+\omega_{i,n+1}N,
\qquad \omega_{i,n+1}=\sum_jh_{ij}\omega_j.
\]
Thus $h_{ij}=h_{ji}$, $H=\sum_i h_{ii}$, and
\[
\vec h=\sum_{i,j}h_{ij}\omega_i\otimes\omega_jN.
\]
Here $\vec h$ is the normal-valued second fundamental form. We also
write $h$ for its scalar form and the associated shape operator,
$h(e_i)=\sum_jh_{ij}e_j$. Thus $\vec h=hN$. The Weingarten formula is
$dN=-dX\circ h$. Matrix products and traces refer to $(h_{ij})$.
For a function $v$, we write $v_{,i}=\nabla_i v$ and
$v_{,ij}=(\nabla^2v)(e_i,e_j)$.
We denote the Euclidean connection by $\overline\nabla$.
Let $\lambda_1,\ldots,\lambda_n$ denote the principal curvatures,
that is, the eigenvalues of $(h_{ij})$. Following \cite{CLW1}, define
\begin{equation}\label{eq:curvature-invariants}
\begin{split}
S&=\sum_{i,j}h_{ij}^2=\sum_i\lambda_i^2, \ \ 
f_3=\sum_{i,j,k}h_{ij}h_{jk}h_{ki}=\sum_i\lambda_i^3,\ \
f_4=\sum_{i,j,k,l}h_{ij}h_{jk}h_{kl}h_{li}=\sum_i\lambda_i^4.
\end{split}
\end{equation}
The Gauss equations are
\begin{equation}\label{eq:gauss}
\begin{split}
R_{ijkl}&=h_{ik}h_{jl}-h_{il}h_{jk},\ \ 
R_{ij}=Hh_{ij}-\sum_kh_{ik}h_{kj},\quad R=H^2-S.
\end{split}
\end{equation}
Here $R_{ij}$ and $R$ denote the Ricci tensor and scalar curvature,
respectively.
\\
We use the derivative order of \cite{CLW1}:
\begin{align}
\sum_kh_{ijk}\omega_k
 &=dh_{ij}+\sum_kh_{ik}\omega_{kj}+\sum_kh_{kj}\omega_{ki},\label{eq:derivone}\\
\sum_lh_{ijkl}\omega_l
 &=dh_{ijk}+\sum_lh_{ljk}\omega_{li}
   +\sum_lh_{ilk}\omega_{lj}+\sum_lh_{ijl}\omega_{lk}.\label{eq:derivtwo}
\end{align}
The Codazzi equations imply that $h_{ijk}$ is symmetric in all three
indices. The Ricci identity is
\begin{equation}\label{ricci-commutation}
h_{ijkl}-h_{ijlk}
=\sum_mh_{mj}R_{mikl}+\sum_mh_{im}R_{mjkl}.
\end{equation}
Following \cite[Section~3]{CLW1}, define
\begin{equation}\label{eq:AB}
\begin{split}
A&=\sum_{i,j,k,l,m}h_{ijk}h_{ljk}h_{im}h_{ml}=\sum_{i,j,k}\lambda_i^2h_{ijk}^2,\ \
B=\sum_{i,j,k,l,m}h_{ijk}h_{lmk}h_{il}h_{jm}=\sum_{i,j,k}\lambda_i\lambda_jh_{ijk}^2.
\end{split}
\end{equation}
Then
\begin{equation}\label{eq:Ddef}
A+2B
=\frac13\sum_{i,j,k}(\lambda_i+\lambda_j+\lambda_k)^2h_{ijk}^2\ge0.
\end{equation}
\\
We also use 
\begin{equation}\label{eq:CD}
\begin{split}
 C&=\sum_{i,j,k,l}h_{ij}h_{jlk}h_{lik}
   =\sum_{i,j,k}\lambda_i h_{ijk}^2,\ \
 D_k=\sum_{i,j,l}h_{ij}h_{jl}h_{lik}
   =\sum_i\lambda_i^2h_{iik}.
\end{split}
\end{equation}
It is easy to get 
\begin{equation}\label{eq:trace-gradients}
 (f_3)_{,k}=3D_k,\qquad (f_4)_{,k}=4\sum_i\lambda_i^3h_{iik}.
\end{equation}
For $S>0$, define
\begin{equation}\label{eq:f-def}
 f=f_4-\frac{f_3^2}{S}.
\end{equation}
We recall \cite[equation~(3.1)]{CWY}:
\begin{equation}\label{fidentity}
  Sf=\sum_{i<j}\lambda_i^2\lambda_j^2(\lambda_i-\lambda_j)^2,
 \qquad f\geq 0.
\end{equation}
\\
We put
\begin{equation}\label{eq:drift-definition}
 \psi=\frac{|X|^2}{2},\qquad \rho=e^{-\psi},\qquad
 \Lop v=\Delta v-\langle X, \nabla v\rangle
\end{equation}
for a smooth function $v$. For a tangent vector field $Z$, define
\begin{equation}\label{eq:weighted-divergence}
 \Div_\rho Z=\rho^{-1}\Div(\rho Z)
 =\Div Z-\langle\nabla\psi,Z\rangle.
\end{equation}
In particular, $\Div_\rho(\nabla v)=\Lop v$.
If $\eta$ has compact support, the divergence theorem gives
\begin{equation}\label{eq:weighted-parts}
 \int_M\eta\Div_\rho Z\,\rho\,d\mu
 =-\int_M\langle\nabla\eta,Z\rangle\rho\,d\mu.
\end{equation}
This identity also holds for compactly supported Lipschitz $\eta$ by
smooth approximation on its compact support.
\\
Put $x_i=\langle X,e_i\rangle$. Differentiating \eqref{shrinker}
and using $\overline\nabla_{e_i}N=-\sum_jh_{ij}e_j$, we obtain
\begin{equation}\label{eq:first-position}
 H_{,i}=\sum_jh_{ij}x_j,\qquad (x_i)_{,j}=\delta_{ij}-Hh_{ij}.
\end{equation}
Consequently,
\begin{equation}\label{eq:psiderivs}
 \nabla\psi=X^\top,\qquad \psi_{,ij}=\delta_{ij}-Hh_{ij},
 \qquad \nabla H=hX^\top,
\end{equation}
where  $X^\top$ denotes the tangential component of $X$.
A second differentiation of \eqref{eq:first-position}, followed by
Codazzi, gives
\begin{equation}\label{mean-curvature-hessian}
 H_{,ij}=\sum_k h_{ijk}x_k+h_{ij}-H\sum_kh_{ik}h_{kj}.
\end{equation}
\\
The following standard identities will be used throughout the paper;
see Cheng--Wei \cite{CW}, Cheng--Wei--Yano
\cite[equations~(2.5)--(2.6) and Lemma~2.1]{CWY}, and
Cheng--Li--Wei \cite[Section~3]{CLW1}.
\begin{lemma}[\cite{CW,CWY,CLW1}]\label{basic-identities}
For a smooth self-shrinker, we have
\begin{equation}\label{eq:Lh}
 \Lop h_{ij}=(1-S)h_{ij},\qquad \Lop H=(1-S)H,
\end{equation}
\begin{equation}\label{eq:LS}
 \frac12\Lop S=\sum_{i,j,k}h_{ijk}^2+S(1-S),
\end{equation}
\begin{equation}\label{evolutions}
 \frac13\Lop f_3=(1-S)f_3+2C,\qquad
 \frac14\Lop f_4=(1-S)f_4+2A+B,
\end{equation}
\begin{equation}\label{eq:LHsquare}
 \frac12\Lop H^2=|\nabla H|^2+(1-S)H^2,
\end{equation}
and
\begin{equation}\label{eq:Lposition}
 \Lop X=-X,\qquad \frac12\Lop|X|^2=n-|X|^2,
 \qquad \Lop\psi=n-2\psi.
\end{equation}
If $S$ is constant, then
\begin{equation}\label{simons}
 \sum_{i,j,k}h_{ijk}^2=S(S-1)
\end{equation}
and
\begin{equation}\label{secondderivative}
 \sum_{i,j,k,l}h_{ijkl}^2=S(S-1)(S-2)+3(A-2B).
\end{equation}
\end{lemma}
It is obvious that if $S$ is constant and $S\le1$, then $S=0$ or $S=1$, and $\nabla h=0$
from  \eqref{simons}.
\\
When $S$ is constant, Cauchy--Schwarz and \eqref{simons} give
\begin{equation}\label{eq:Hbounds}
 H^2\le nS,\qquad
 |\nabla H|^2=\sum_k\left(\sum_i h_{iik}\right)^2
 \le n\sum_{i,j,k}h_{ijk}^2=nS(S-1).
\end{equation}
Furthermore,
\begin{equation}\label{eq:ABbounds}
 0\le A\le S\sum_{i,j,k}h_{ijk}^2,\qquad |B|\le S\sum_{i,j,k}h_{ijk}^2.
\end{equation}
It follows from \eqref{simons}, \eqref{secondderivative} and
\eqref{eq:ABbounds} that $h$, $\nabla h$ and $\nabla^2h$ are uniformly
bounded. For the estimates in which $S>1$, put
\begin{equation}\label{eq:t-definition}
 t=\frac{S-1}{S},\qquad \sum_{i,j,k}h_{ijk}^2=tS^2.
\end{equation}
\\
The standard mean curvature Hessian formula in
\cite[equations~(3.5)--(3.6)]{CP}
gives the following identities.
\begin{lemma}[\cite{CP}]\label{mean-curvature-contractions}
If $S$ is constant, then, in a principal frame,
\begin{align}\label{mean-curvature-moment-identities}
 \sum_i\lambda_iH_{,ii}&=S-Hf_3,\notag\\
 \sum_i\lambda_i^2H_{,ii}
 &=f_3-Hf_4+\frac13\langle X,\nabla f_3\rangle.
\end{align}
\end{lemma}

For constant $S>1$, we use the derivative estimate of
Cheng--Wei \cite[Lemma~3.5]{CW}; see also
\cite[equation~(3.1)]{CLW1}. With the notation above, it reads as follows.
\begin{lemma}[{\cite[Lemma~3.5]{CW}}]\label{bessel}
If $S>1$ is constant, then, in a principal frame, we have
\begin{equation}\label{besselbounds}
3C^2\le tS^2
 \left(A+2B-4\sum_k\frac{D_k^2}{S+2\lambda_k^2}\right).
\end{equation}
\end{lemma}

If $S>1$ and $f_3$ are constant, Cheng--Wei--Yano
\cite[equations~(3.2)--(3.3)]{CWY} give
\begin{equation}\label{c3:constant-norm-C}
 2C=(S-1)f_3=tSf_3,
\end{equation}
\begin{equation}\label{c3:ABidentity}
 A-2B=Sf=Sf_4-f_3^2.
\end{equation}
Their differentiated constant trace identities give
\begin{equation}\label{eq:constant-trace-gradients}
 \sum_i\lambda_i h_{iik}=0,\qquad
 \sum_i\lambda_i^2h_{iik}=D_k=0,
\end{equation}
\begin{align}
 \sum_i\lambda_i h_{iikl}&=-\sum_{i,j}h_{ijk}h_{ijl},\label{c3:HessS}\\
 \sum_i\lambda_i^2h_{iikl}&=-2\sum_{i,j}\lambda_i h_{ijk}h_{ijl}.
 \label{c3:Hessf3}
\end{align}
The Ricci identity also gives \cite[equation~(2.9)]{CWY}
\begin{equation}\label{eq:principal-ricci}
 h_{iijj}-h_{jjii}=\lambda_i\lambda_j(\lambda_i-\lambda_j).
\end{equation}

\begin{proposition}
\label{cubic-sup-estimate}
Let $X:M^n\to\R^{n+1}$ be a  complete self-shrinker.
Suppose that $S,f_3$ are constant and $1<S<4$. Then $f_3\ne0$ and
\begin{equation}\label{cubic-sup-estimate-bound}
 S\sup_M f\le\frac{S-1}{4(4-S)}f_3^2.
\end{equation}
\end{proposition}
\begin{proof}
Since $f_3$ is constant, $\vec D=0$. \eqref{besselbounds} and \eqref{c3:constant-norm-C} yield
\[
 A+2B\ge\frac{3t}{4}f_3^2.
\]
Since $S$ is constant, Gauss equation \eqref{eq:gauss} yields that Ricci curvature is bounded and
$f_4$, $h_{ijk}$ and $h_{ijkl}$ for any $i, j, k, l$ are bounded. By applying the generalized maximum principle 
of Cheng and Peng \cite{CP} recalled in Subsection~\ref{subsec:maxprinciples} to $f_4$, we get that
there is a sequence $p_m\in M$ such that
\begin{equation*}
f_4(p_m)\to\sup_M f_4,\qquad |\nabla f_4|(p_m)\to0,
\qquad \limsup_{m\to\infty}\Lop f_4(p_m)\le0.
\end{equation*}
Since 
$f_4$, $h_{ijk}$ and $h_{ijkl}$ for any $i, j, k, l$ are bounded, we can assume 
that $\lambda_i(p_m)$, $h_{ijk}(p_m)$ and $h_{ijkl}(p_m)$ converge if necessary taking a subsequence.
Hence, we know 
$f_4\to \sup_M f_4$ and $A,B$ have limits $A_*,B_*$. \eqref{evolutions} yields 
\begin{equation}\label{c3:maximum-ineq}
 2A_*+B_*\le tS\sup_M f_4.
\end{equation}
In view of  $\sup_M f=\sup_M f_4-f_3^2/S$ and \eqref{c3:ABidentity}, we infer
\begin{align*}
 4t(S\sup_M f+f_3^2)&\ge4(2A_*+B_*)\\
 &=3(A_*-2B_*)+5(A_*+2B_*)\\
 &\ge3S\sup_M f+\frac{15t}{4}f_3^2.
\end{align*}
Hence
\[
 (3-4t)S\sup_M f\le\frac t4 f_3^2,
 \qquad 3-4t=\frac{4-S}{S}>0.
\]
This proves \eqref{cubic-sup-estimate-bound}. If $f_3=0$, then \eqref{cubic-sup-estimate-bound} yields
\[
 0\le f_4=f\le\sup_M f=0,\qquad S=0,
\]
contrary to $S>1$.
\end{proof}

\subsection{Maximum principles}\label{subsec:maxprinciples}
We will use the following maximum principles.

\begin{theorem}[Cheng--Peng generalized maximum principle
{\cite[Theorem~3.1]{CP}}]\label{thm:CP}
Let $X:M^n\to\R^{n+1}$ be a complete self-shrinker whose Ricci curvature
is bounded below. For any $v\in C^2(M)$ bounded above, there is a sequence
$p_m\in M$ such that
\begin{equation}\label{eq:CP}
v(p_m)\to\sup_M v,\qquad |\nabla v|(p_m)\to0,
\qquad \limsup_{m\to\infty}\Lop v(p_m)\le0.
\end{equation}
\end{theorem}

\begin{theorem}[Hessian Omori--Yau principle {\cite{Omori}}]\label{thm:Hessian}
Let $(M,g)$ be complete with sectional curvature bounded below.
For any $v\in C^2(M)$ bounded below, there are points $p_m\in M$ such that
\begin{equation}\label{eq:Hessianprinciple}
\begin{split}
&v(p_m)\to\inf_M v,\qquad |\nabla v|(p_m)\to0,\\
&\nabla^2v(p_m)(w,w)\ge-m^{-1}|w|^2
\qquad\text{for every }w\in T_{p_m}M.
\end{split}
\end{equation}
For a bounded-above function, one obtains a supremum sequence with
$\nabla^2v(p_m)\le m^{-1}g$.
\end{theorem}

\begin{remark}
\label{rem:maxprinciples}
For complete self-shrinkers with constant $S$, \eqref{eq:gauss} and \eqref{eq:Hbounds} give
\begin{equation}\label{eq:curvature-lower-bounds}
\operatorname{Ric}\ge-(\sqrt n+1)Sg,\qquad
\sec(e_i, e_j)=\lambda_i\lambda_j\ge-\frac{\lambda_i^2+\lambda_j^2}{2}\ge-\frac S2.
\end{equation}
Theorem~\ref{thm:CP} gives an estimate for $\Lop v$.
Theorem~\ref{thm:Hessian} gives an estimate for $\nabla^2v$.
If $X$ is unbounded, the latter does not by itself control
$\langle X,\nabla v\rangle$. We use these two principles separately.
\end{remark}

We also use the distance-barrier form of the Omori-Yau principle.
The trace-operator conclusion below is a special case of the exhaustion
criterion of Albanese--Al\'ias--Rigoli
\cite[Theorems~A--B and Remark~3.1]{AAR}.
We record the simultaneous Hessian and gradient estimates needed later.

\begin{lemma}[Distance-barrier consequence of \cite{Omori,AAR}]
\label{trace-minimum-principle}
Let $M$ be complete with sectional curvature bounded below.
Fix $o\in M$ and put $r=d(o,\cdot)$.
For every bounded below $v\in C^2(M)$, there are $p_j\in M$ and
$\varepsilon_j\downarrow0$ such that
\begin{equation}\label{eq:simultaneous-minimum}
\begin{split}
 v(p_j)&\longrightarrow\inf_M v,\\
 \nabla^2v(p_j)&\ge-\varepsilon_j g,\qquad
 (1+r(p_j))|\nabla v|(p_j)\le\varepsilon_j.
\end{split}
\end{equation}
In particular, if $T\ge0$ is symmetric with bounded trace and
$|W|\le K(1+r)$, then
\[
 Lv=\tr(T\circ\nabla^2v)+\langle W,\nabla v\rangle
 \quad\Longrightarrow\quad \liminf_j Lv(p_j)\ge0.
\]
\end{lemma}
\begin{proof}[Verification of the distance-barrier estimates]
For compact $M$, take a minimum point. Otherwise set
\[
 \Phi=\log(1+r^2),\qquad \operatorname{Sec}\ge-k^2.
\]
Away from the cut locus, Hessian comparison gives
\begin{align*}
 |\nabla\Phi|&=\frac{2r}{1+r^2}\le\frac4{1+r},\\
 \nabla^2\Phi
 &=\frac{2r}{1+r^2}\nabla^2r
   +\frac{2(1-r^2)}{(1+r^2)^2}dr\otimes dr
 \le(k+4)g
\end{align*}
because of $\nabla^2r\leq (k+\dfrac1r)g$.
At a cut point, use the Calabi upper support
\cite[Remark~3.1]{AAR}:
\[
 r_\delta(q)=\delta+d(o_\delta,q),\qquad 0<\delta<r/2,
 \qquad \Phi_\delta=\log(1+r_\delta^2).
\]
At contact, $r_\delta=r$ and
\[
 |\nabla\Phi_\delta|=\frac{2r}{1+r^2},\qquad
 \nabla^2\Phi_\delta
 \le\left\{\frac{2r}{1+r^2}
       \left(k+\frac1{r-\delta}\right)+2\right\}g
 \le(k+6)g.
\]
The function $\Phi$ is smooth near $o$ and proper.
Let $p_\eta$ minimize $v+\eta\Phi$. Using the upper supports if needed,
\[
 \nabla^2v(p_\eta)\ge-(k+6)\eta g,\qquad
 (1+r(p_\eta))|\nabla v|(p_\eta)\le4\eta.
\]
For every fixed $q\in M$,
\[
 \inf_Mv\le v(p_\eta)\le v(q)+\eta\Phi(q).
\]
Let $\eta\downarrow0$ and then $v(q)\downarrow\inf_Mv$.
After rescaling $\eta$, this proves \eqref{eq:simultaneous-minimum}.
Finally,
\[
 Lv(p_j)\ge-\varepsilon_j\bigl(\sup_M\tr T+K\bigr)
 \longrightarrow0.
\]
\end{proof}

\subsection{Sub-exponential volume growth}\label{subsec:subexp}

Throughout this subsection, $(M^n,g)$ is connected and complete.
By Hopf--Rinow theorem, every finite-radius intrinsic ball has finite volume.
We use the following condition; see also Brooks \cite{Brooks81} and
Sturm \cite{Sturm93}.

\begin{definition}[\cite{Brooks81,Sturm93}]\label{def:subexp}
We say that $M$ has sub-exponential volume growth if, for some $p\in M$
and every $\epsilon>0$, there is a constant $C_\epsilon>0$ such that
\[
\operatorname{Vol}(B_r^M(p))\le C_\epsilon e^{\epsilon r}
\qquad\text{for all }r\ge1.
\]
The constant may depend on $\epsilon$ and $p$, but not on $r$.
\end{definition}

The uniform version considered in \cite{Sturm93} requires
\[
\operatorname{Vol}(B_r^M(p))
\le C_\epsilon e^{\epsilon r}\operatorname{Vol}(B_1^M(p))
\qquad(p\in M,\ r\ge1)
\]
with the same constant $C_\epsilon$ for all $p$.
In Definition~\ref{def:subexp}, $C_\epsilon$ may depend on $p$.

\begin{proposition}
\label{prop:subexp-equivalence}
For a fixed $p\in M$, sub-exponential volume growth is equivalent to
\begin{equation}\label{eq:subexp-limit}
\limsup_{r\to\infty}\frac{\log\operatorname{Vol}(B_r^M(p))}{r}=0.
\end{equation}
In this case the limit exists and equals zero. The condition is independent
of the base point.
\end{proposition}
\begin{proof}
The defining bound gives
\[
\frac{\log\operatorname{Vol}(B_r^M(p))}{r}
\le\epsilon+\frac{\log C_\epsilon}{r}.
\]
Hence the upper limit is at most zero. By monotonicity,
\[
\operatorname{Vol}(B_r^M(p))\ge\operatorname{Vol}(B_1^M(p))>0
\qquad(r\ge1).
\]
Thus the lower limit is at least zero.

Now assume \eqref{eq:subexp-limit}. Given $\epsilon>0$, choose $R\ge1$
so that
\[
\operatorname{Vol}(B_r^M(p))\le e^{\epsilon r}\quad(r\ge R).
\]
For $1\le r\le R$, use
$\operatorname{Vol}(B_r^M(p))\le\operatorname{Vol}(B_R^M(p))$.
Thus Definition~\ref{def:subexp} holds with
$C_\epsilon=\max\{1,\operatorname{Vol}(B_R^M(p))\}$.

For $q\in M$, put $d=d_M(p,q)$. The triangle inequality yields
\begin{equation}\label{eq:ball-shift}
B_{r-d}^M(p)\subset B_r^M(q)\subset B_{r+d}^M(p)
\qquad(r>d).
\end{equation}
The second inclusion holds for every $r>0$. It gives
\[
\operatorname{Vol}(B_r^M(q))
\le C_\epsilon e^{\epsilon(r+d)}
=(C_\epsilon e^{\epsilon d})e^{\epsilon r}.
\]
Interchanging $p$ and $q$, we prove the assertion.
\end{proof}

\begin{remark}\label{rem:subexp-polynomial}
At most Euclidean volume growth means
$\operatorname{Vol}(B_r^M(p))\le C(1+r)^n$. Polynomial volume growth allows
any fixed exponent $d\ge0$. Hence
\[
\begin{aligned}
&\text{at most Euclidean volume growth}\ \Longrightarrow\
\text{polynomial volume growth}\\
& \Longrightarrow\ 
\text{sub-exponential volume growth}.
\end{aligned}
\]
Indeed,
\[
\frac{\log\operatorname{Vol}(B_r^M(p))}{r}
\le\frac{\log C+d\log(1+r)}r\longrightarrow0
\]
and Proposition~\ref{prop:subexp-equivalence} applies. The converse
need not hold for growth functions: $e^{\sqrt r}$ is sub-exponential, but
\[
\sqrt r-d\log r\longrightarrow\infty\qquad(d\ge0).
\]
Thus a sub-exponential growth function need not be polynomial.
An estimate with one fixed positive exponential rate is also
insufficient for Definition~\ref{def:subexp}.
\end{remark}

The following comparison distinguishes intrinsic volume growth from
extrinsic volume growth.
\begin{proposition}\label{prop:subexp-extrinsic}
Let $X:M^n\to\R^{n+1}$ be an isometric immersion. Any polynomial volume 
bound on $\operatorname{Vol}(X^{-1}(B_R^{\R^{n+1}}(0)))$ implies an
intrinsic polynomial volume bound of the same degree. The same implication holds
for sub-exponential volume growth. Neither converse holds for general immersions.
\end{proposition}
\begin{proof}
Fix $p\in M$. For any piecewise smooth curve $\gamma$ joining $p$
to $q$, the isometric immersion satisfies
\[
|X(q)-X(p)|\le L_{\R^{n+1}}(X\circ\gamma)=L_M(\gamma).
\]
Taking the infimum over $\gamma$, we obtain
\[
|X(q)|\le |X(p)|+d_M(p,q).
\]
It follows that
\[
B_r^M(p)\subset
X^{-1}\bigl(B_{r+|X(p)|}^{\R^{n+1}}(0)\bigr).
\]
Suppose first that the extrinsic volume satisfies
\[
\operatorname{Vol}\bigl(X^{-1}(B_R^{\R^{n+1}}(0))\bigr)
\le C(1+R)^d,\qquad R\ge1.
\]
For $r\ge1$, we then have
\begin{align*}
\operatorname{Vol}(B_r^M(p))
&\le C(1+r+|X(p)|)^d\\
&\le C(1+|X(p)|)^d(1+r)^d.
\end{align*}
This proves the polynomial bound with the same exponent.
If the extrinsic volume has sub-exponential growth, then for every
$\epsilon>0$,
\begin{align*}
\operatorname{Vol}(B_r^M(p))
&\le C_\epsilon e^{\epsilon(r+|X(p)|)}\\
&=\bigl(C_\epsilon e^{\epsilon|X(p)|}\bigr)e^{\epsilon r}.
\end{align*}
Thus the intrinsic volume also has sub-exponential growth.

For the converse, consider the example in Example~\ref{rem:cover}:
\[
X(t,z)=(\cos t,\sin t,z),\qquad
(t,z)\in\R\times\R^{n-1}.
\]
Its induced metric is $dt^2+|dz|^2$. Hence
\[
\operatorname{Vol}(B_r^M(p))=\omega_n r^n.
\]
On the other hand, $|X(t,z)|^2=1+|z|^2$. For $R>1$,
\[
X^{-1}(B_R^{\R^{n+1}}(0))
=\R\times B_{\sqrt{R^2-1}}^{\R^{n-1}}(0)
\]
has infinite volume. This excludes both converse implications.
\end{proof}

For weighted sub-exponential growth and exponential kernel estimates,
see \cite[Definition~6.2 and Lemma~6.3]{CharLu15}.
Their result assumes uniform growth in the center. The following
statement only requires growth at a fixed base point. We give its
proof under this assumption.

\begin{proposition}
\label{prop:subexp-weighted}
Let $d\nu=w\,d\mu$, where $w$ is smooth and positive. Then
Proposition~\ref{prop:subexp-equivalence} holds with $\nu(B_r^M(p))$ in
place of ordinary volume. Moreover, sub-exponential volume growth of these
weighted balls is equivalent to
\begin{equation}\label{eq:subexp-integrability}
\int_M e^{-a d_M(p,x)}\,d\nu(x)<\infty
\qquad\text{for every }a>0.
\end{equation}
For an immersion, ordinary sub-exponential volume growth implies Gaussian
sub-exponential volume growth, with $w=\rho=e^{-|X|^2/2}$.
\end{proposition}
\begin{proof}
Put
\[
r(x)=d_M(p,x),\qquad W_p(R)=\nu(B_R^M(p)).
\]
By completeness, $\overline{B_R^M(p)}$ is compact. Since $w$ is smooth
and positive,
\[
0<W_p(R)<\infty\quad(R>0),\qquad
W_p(R)\ge W_p(1)>0\quad(R\ge1).
\]
If $W_p(R)\le C_\epsilon e^{\epsilon R}$ for every $\epsilon>0$, then
\[
\frac{\log W_p(1)}{R}
\le\frac{\log W_p(R)}{R}
\le\epsilon+\frac{\log C_\epsilon}{R},\qquad R\ge1.
\]
Letting $R\to\infty$ and then $\epsilon\downarrow0$, we obtain
$\log W_p(R)/R\to0$. Conversely, this limit gives
$W_p(R)\le e^{\epsilon R}$ for all sufficiently large $R$.
Increasing the constant covers $1\le R$ below that radius.
Thus the weighted version of \eqref{eq:subexp-limit} holds.
For $q\in M$, \eqref{eq:ball-shift} gives
\[
W_q(R)\le W_p(R+d_M(p,q))
\le C_\epsilon e^{\epsilon d_M(p,q)}e^{\epsilon R}.
\]
Interchanging $p$ and $q$ proves independence of the base point.
\\
For every $a>0$, we have
\[
e^{-ar(x)}=\int_{r(x)}^\infty ae^{-at}\,dt.
\]
The integrand is nonnegative. By Tonelli's theorem,
\begin{align*}
\int_M e^{-ar(x)}\,d\nu(x)
&=\int_M\int_{r(x)}^\infty ae^{-at}\,dt\,d\nu(x)\\
&=a\int_0^\infty e^{-at}W_p(t)\,dt.
\end{align*}
If $W_p$ has sub-exponential growth, choose $\epsilon=a/2$. Then
\begin{align*}
\int_M e^{-ar(x)}\,d\nu(x)
&\le aW_p(1)\int_0^1e^{-at}\,dt
  +aC_{a/2}\int_1^\infty e^{-at/2}\,dt\\
&\le W_p(1)+2C_{a/2}e^{-a/2}<\infty.
\end{align*}
This proves \eqref{eq:subexp-integrability}.
Conversely, suppose that the integral in
\eqref{eq:subexp-integrability} is finite for every $a>0$.
Since $e^{-ar(x)}\ge e^{-aR}$ on $B_R^M(p)$,
\[
e^{-aR}W_p(R)
\le\int_{B_R^M(p)}e^{-ar(x)}\,d\nu(x)
\le\int_Me^{-ar(x)}\,d\nu(x).
\]
Multiplication by $e^{aR}$ gives the sub-exponential bound with
$\epsilon=a$. This proves the equivalence.

For the Gaussian measure, $0<\rho\le1$ implies
\[
\int_{B_R^M(p)}\rho\,d\mu
\le\operatorname{Vol}(B_R^M(p)).
\]
Hence ordinary sub-exponential volume growth implies Gaussian weighted
sub-exponential volume growth. This completes the proof.
\end{proof}

\begin{remark}
\label{rem:weightedgrowthscope}
The infinite cylinder covering in Example~\ref{rem:cover} is not
proper. Its intrinsic balls have polynomial volume growth, and hence
sub-exponential volume growth. By Proposition~\ref{prop:subexp-weighted},
its Gaussian weighted intrinsic balls also have sub-exponential volume growth.
Thus neither growth condition in Theorem~\ref{thm:maingrowth}
implies properness of the immersion.
\end{remark}

\section{Proof of Theorem~\ref{thm:f3}}\label{sec:algebra-cubic}
We assume $S>1$ throughout this section.
Constancy of $f_3$ is imposed where stated.

\begin{proposition}\label{mean-curvature-changes-sign}
If $S$ is constant and $S>1$, then
\[
 \inf_M H<0<\sup_M H.
\]
Only constant $S$ is required in this proposition.
\end{proposition}
\begin{proof}
Suppose that $H\ge0$. 
By Lemma~\ref{trace-minimum-principle}, there exist $p_m\in M$ and
$\varepsilon_m\downarrow0$ such that
\[
 H(p_m)\to \inf_M H,\qquad \nabla^2H(p_m)\ge-\varepsilon_m g,
 \qquad |\langle X,\nabla H\rangle|(p_m)\le\varepsilon_m.
\]
Here we used $\operatorname{Sec}\ge-S/2$ and
$|X(p)|\le|X(o)|+d(o,p)$.
Choose a principal frame at $p_m$. Put
\[
 a_i^{(m)}=H_{,ii}(p_m)+\varepsilon_m\ge0.
\]
According to \eqref{eq:Lh},
\[
 \sum_i a_i^{(m)}
 =(1-S)H(p_m)+\langle X,\nabla H\rangle(p_m)+n\varepsilon_m
 \longrightarrow(1-S)\inf_M H.
\]
Thus $(1-S)\inf_M H\ge0$. Since $S>1$ and $\inf_M H\ge0$, we have
\[
 \inf_M H=0,\qquad \sum_i a_i^{(m)}\longrightarrow0.
\]
Since $|\lambda_i|\le\sqrt S$,
\[
 \left|\sum_i\lambda_i(p_m)a_i^{(m)}\right|
 \le\sqrt S\sum_i a_i^{(m)}\longrightarrow0.
\]
On the other hand, \eqref{mean-curvature-moment-identities} gives
\[
 \sum_i\lambda_i(p_m)a_i^{(m)}
 =S-H(p_m)f_3(p_m)+\varepsilon_mH(p_m)\longrightarrow S,
\]
where $|f_3|\le S^{3/2}$. This is a contradiction.
Hence $\inf_M H<0$. Reversing the normal, we obtain $\sup_M H>0$.
\end{proof}

\begin{proposition}\label{mean-curvature-discriminant}
Suppose that $S,f_3$ are constant and $S>1$. Then
\begin{equation}\label{negative-discriminant-necessary}
 (S-1)\sup_M f_4
 \ge f_3^2+\frac{S+1}{\sqrt{nS}}|f_3|+\frac Sn.
\end{equation}
We also have
\begin{equation}\label{constant-cubic-amplitude-bound}
 |f_3|\le\sqrt S\left(S-1-\frac1{\sqrt n}\right)
\end{equation}
and \begin{equation}\label{spectral-cubic-amplitude}
 |f_3|\le\sqrt S\left(S-1-2\sqrt{\frac{2(S-1)}n}\right).
\end{equation}
In particular, $(S-1)\sup_M f_4>f_3^2$.
\end{proposition}
\begin{proof}
Since $f_3$ is constant, we can assume  $f_3\ge0$. 
Proposition~\ref{mean-curvature-changes-sign} and \eqref{eq:Hbounds} give
\[
 0<-\inf_M H\le\sqrt{nS}.
\]
By Lemma~\ref{trace-minimum-principle}, there exist $p_m\in M$ and
$\varepsilon_m\downarrow0$ such that
\[
 H(p_m)\to\inf_M H,\qquad \nabla^2H(p_m)\ge-\varepsilon_mg,
 \qquad |\langle X,\nabla H\rangle|(p_m)\le\varepsilon_m.
\]
Choose a principal frame at $p_m$, and put
$a_i^{(m)}=H_{,ii}(p_m)+\varepsilon_m\ge0$.
By \eqref{eq:Lh} and \eqref{mean-curvature-moment-identities},
\begin{equation}\label{eq:3.0}
\begin{aligned}
 \sum_i a_i^{(m)}
 &=(1-S)H+\langle X,\nabla H\rangle+n\varepsilon_m,\\
 \sum_i\lambda_i a_i^{(m)}&=S-Hf_3+\varepsilon_mH,\\
 \sum_i\lambda_i^2 a_i^{(m)}&=f_3-Hf_4+\varepsilon_mS.
\end{aligned}
\end{equation}
All quantities on the right are evaluated at $p_m$.
Since $\lambda_i\le\sqrt S$,
\[
 \sum_i\lambda_i a_i^{(m)}\le\sqrt S\sum_i a_i^{(m)}.
\]
Taking the limit, we obtain
\[
 S-\inf_M H f_3\le-\sqrt S(S-1)\inf_M H.
\]
Consequently,
\[
 f_3\le\sqrt S(S-1)-\frac S{-\inf_M H}
 \le\sqrt S\left(S-1-\frac1{\sqrt n}\right).
\]
This proves \eqref{constant-cubic-amplitude-bound}.
\\
Cauchy--Schwarz  inequality gives, from $a_i^{(m)}\geq 0$, 
\[
 \left(\sum_i\lambda_i a_i^{(m)}\right)^2
 \le\left(\sum_i a_i^{(m)}\right)
       \left(\sum_i\lambda_i^2a_i^{(m)}\right).
\]
Since $f_4\le\sup_M f_4<\infty$, taking the upper limit in \eqref{eq:3.0} yields
\[
 (S-\inf_M H f_3)^2\le-(S-1)\inf_M H\{f_3-\inf_M H\sup_M f_4\}.
\]
Expanding and dividing by $(\inf_M H)^2$, we have
\begin{align*}
 (S-1)\sup_M f_4-f_3^2
 &\ge\frac{S+1}{-\inf_M H}f_3+\frac{S^2}{(\inf_M H)^2}\\
 &\ge\frac{S+1}{\sqrt{nS}}f_3+\frac Sn.
\end{align*}
This proves \eqref{negative-discriminant-necessary}.
\\
Under sense of limits, let $\sum_i\lambda_i=\inf_M H<0$ and $\sum_i\lambda_i^2=S$.
Put $\lambda_{\max}=\max_i\lambda_i$.
For $n\ge2$, Cauchy--Schwarz  inequality gives
\[
 (\lambda_{\max}-\inf_M H)^2\le(n-1)(S-\lambda_{\max}^2).
\]
Hence
\[
 \lambda_{\max}\le\frac{\inf_M H+\sqrt{(n-1)(nS-(\inf_M H)^2)}}n.
\]
Since  $0\le \dfrac{-\inf_M H}{\sqrt S}\le\sqrt n\le n$, we know 
\[
 (n+\dfrac{\inf_M H}{\sqrt S})^2-(n-1)\biggl\{n-\bigl(\dfrac{\inf_M H}{\sqrt S}\bigl)^2\biggl\}=n(\dfrac{\inf_M H}{\sqrt S}+1)^2\ge0.
\]
Thus, we have
\[
 (n+\dfrac{\inf_M H}{\sqrt S})^2\geq (n-1)\biggl\{n-\bigl(\dfrac{\inf_M H}{\sqrt S}\bigl)^2\biggl\}.
\]
We obtain
\begin{equation}\label{negative-trace-largest-root}
 \lambda_{\max}\le\sqrt S+\frac{2\inf_M H}n.
\end{equation}
For $n=1$, this is an equality.
\\
Since $\{\lambda_i(p_m)\}$ and $\{a^m_i\}$ converge, 
taking  the limit for  the second formula in  \eqref{eq:3.0} 
and 
\eqref{negative-trace-largest-root}  give
\[
 S-\inf_M H f_3=\sum_i\lim\lambda_i (p_m)\lim a^m_i\le\left(\sqrt S+\frac{2\inf_M H}n\right)(1-S)\inf_M H.
\]
It follows that
\begin{align*}
 f_3&\le\sqrt S(S-1)+\frac{2(S-1)\inf_M H}n+\frac S{\inf_M H}\\
 &\le\sqrt S(S-1)-2\sqrt{\frac{2S(S-1)}n}.
\end{align*}
Here we used $a^2+b^2\geq 2ab$. 
This proves \eqref{spectral-cubic-amplitude}.
\end{proof}

\begin{lemma}\label{cubic-gap-three}
If $S,f_3$ are constant and $S>1$, then $S>3$.
\end{lemma}
\begin{proof}
Suppose that $1<S\le3$.
The general maximum estimate in
Proposition~\ref{cubic-sup-estimate} gives $f_3\ne0$ and
\[
 S\sup_M f\le\frac{S-1}{4(4-S)}f_3^2.
\]
Since $f_4=f+f_3^2/S$, this implies, everywhere on $M$,
\[
 f_4\le\frac{3(5-S)}{4S(4-S)}f_3^2.
\]
Consequently,
\begin{align*}
 f_3^2-(S-1)f_4
 &\ge\left(1-\frac{3(S-1)(5-S)}{4S(4-S)}\right)f_3^2\\
 &=\frac{(3-S)(S+5)}{4S(4-S)}f_3^2\ge0.
\end{align*}
This contradicts Proposition~\ref{mean-curvature-discriminant}.
\end{proof}

\begin{proof}[Proof of Theorem~\ref{thm:f3}]
Suppose that $S>1$.
For $n=1$, $H^2=S$ makes $H$ constant, contrary to
Proposition~\ref{mean-curvature-changes-sign}.
For $n=2$, the classification of Cheng--Ogata \cite{CO} excludes $S>1$.
For $n=3$, the classification of Cheng--Li--Wei \cite{CLWf3} excludes $S>1$. Hence $n\ge4$.
\\
Lemma~\ref{cubic-gap-three} gives $S>3$.
Put $a=1/\sqrt n\leq \frac12$.
Suppose, contrary to the assertion, that
\[
 3<S\le 3+a<4.
\]
Proposition~\ref{cubic-sup-estimate} gives $f_3\ne0$ and
\[
 \sup_M f_4\le\frac{3(5-S)}{4S(4-S)}f_3^2.
\]
Combining this with \eqref{negative-discriminant-necessary}, we obtain
\[
 \frac{(S+5)(S-3)}{4S(4-S)}f_3^2
 \ge\frac{S+1}{\sqrt{nS}}|f_3|+\frac Sn.
\]
Putting  $z=|f_3|/S^{3/2}>0$,  we get
\begin{equation}\label{normalized-quantitative-cubic}
 \frac{S(S+5)(S-3)}{4(4-S)}
 \ge\frac{(S+1)a}{z}+\frac{a^2}{z^2}.
\end{equation}
By \eqref{spectral-cubic-amplitude} and $S>3$,
\begin{equation}\label{spectral-gap-amplitude-reduction}
 0<z\le\frac{S-1-2a\sqrt{2(S-1)}}S
       \le\frac{S-1-4a}S.
\end{equation}
\\
For $3<S\le s\le 3+a$, from  $S-1-4a>0$, 
\[
 q(s)=\frac{as}{s-1-4a}
\]
is a positive function of $s$ in $[S, 3+a]$ and 
\[
 q'(s)=-\frac{a(1+4a)}{(s-1-4a)^2}<0.
\]
Hence, $q(s)$ is decreasing.
Consequently, in view of \eqref{spectral-gap-amplitude-reduction}, we obtain
\[
 \frac az\ge q(S)\ge q(3+a)
       =\frac{a(3+a)}{2-3a}.
\]
For $3<s<4$,  function
\[
 G(s)=\frac{s(s+5)(s-3)}{4(4-s)}
\]
is an increasing function of $s$ because of  
\[
 G'(s)
 =\dfrac{(3s-5)(s+3)(4-s)+s(s-3)(s+5)}{4(4-s)^2}>0.
\]
Equation \eqref{normalized-quantitative-cubic} and $S>3$ imply
\begin{equation}\label{spectral-gap-rational-comparison}
 \frac{a(3+a)(8+a)}{4(1-a)}=G(3+a)\ge G(S)\geq 4q(3+a)+q(3+a)^2.
\end{equation}
Since $a\leq \frac 12$, we conclude
\begin{equation*}
\begin{aligned}
& 4q(3+a)+q(3+a)^2-\frac{a(3+a)(8+a)}{4(1-a)}\\
&=\frac{4a(3+a)}{2-3a}+\frac{a^2(3+a)^2}{(2-3a)^2}-\frac{a(3+a)(8+a)}{4(1-a)}\\
&=\frac{a^2(3+a)(24-20a-13a^2)}
          {4(1-a)(2-3a)^2}>0,
\end{aligned}
\end{equation*}
because of 
\begin{align*}
 &4(1-a)\{4(2-3a)+a(3+a)\}
                         -(8+a)(2-3a)^2\\
 &\quad=(32-68a+40a^2-4a^3)
                         -(32-92a+60a^2+9a^3)\\
 &\quad=a(24-20a-13a^2)>0.
\end{align*}
This contradicts \eqref{spectral-gap-rational-comparison}.
Hence $S>3+1/\sqrt n$.
\\
If $S\le3+1/\sqrt n$, the preceding conclusion excludes $S>1$.
The formula \eqref{simons}   gives $S=0$ or $S=1$.
Proposition~\ref{prop:parallel-classification} below gives the stated images.
\end{proof}

\begin{proposition}[see \cite{CP,CWY}]\label{prop:parallel-classification}
Let $X:M^n\to\R^{n+1}$ be a connected, complete self-shrinker
with constant $S=0$ or $S=1$. Then, up to an orthogonal transformation,
$X(M)$ is a hyperplane through the origin or
$S^k(\sqrt k)\times\R^{n-k}$ for some $1\le k\le n$.
Moreover, $X$ is a Riemannian covering of its image.
\end{proposition}

\section{Proof of Theorem 1.1}\label{sec:four-dimensional}

\subsection{The signs of the principal curvatures}
We now assume that $S$ and $f_3$ are constant.
No upper bound on $S$ is imposed.

\begin{lemma}[Okumura inequality {\cite{Okumura}}]\label{tracezero-cubic}
Let $n\ge3$ and let $a_1,\ldots,a_n$ be real numbers such that
\[
 \sum_i a_i=0,\qquad \sum_i a_i^2=S>0.
\]
Then
\begin{equation}\label{tracezero-cubic-bound}
 \left(\sum_i a_i^3\right)^2
 \le\frac{(n-2)^2}{n(n-1)}S^3.
\end{equation}
Equality holds precisely when, after
permutation,
\[
 a_1=(n-1)\sqrt{\frac{S}{n(n-1)}},\qquad a_2=\cdots=a_n=-\sqrt{\frac{S}{n(n-1)}}.
 \]
\end{lemma}

\begin{proposition}\label{strict-cubic-obstruction}
Let $n\ge3$. If $S,f_3$ are constant and $S>1$, then
\begin{equation}\label{constant-cubic-upper}
 f_3^2<\frac{(n-2)^2}{n(n-1)}S^3.
\end{equation}
\end{proposition}
\begin{proof}
The assertion is immediate if $f_3=0$.
Without loss of generality, we assume  $f_3>0$.
By Proposition~\ref{mean-curvature-changes-sign}, we have  $\inf_M H<0$.
We consider  the compact set
\[
 \left\{(a_1,\ldots,a_n):\ \sum_i a_i^2=S,\quad
                                      \sum_i a_i\le0\right\}.
\]
We maximize $\sum_i a_i^3$ on this set.
On its boundary $\sum_i a_i=0$, Lemma~\ref{tracezero-cubic} gives
the maximum
\[
 \frac{n-2}{\sqrt{n(n-1)}}S^{3/2}>0.
\]
An interior maximum would satisfy
\[
 3a_i^2=2\alpha a_i.
\]
Thus all nonzero entries would have one common value.
Their sum is negative, so this common value and their cubic sum
would be negative. This cannot be a global maximum.
Consequently, the maximum is attained only where $\sum_i a_i=0$.
At the point of $M$ with $H<0$, we therefore have the strict inequality
\[
 0<f_3<\frac{n-2}{\sqrt{n(n-1)}}S^{3/2}.
\]
Squaring proves \eqref{constant-cubic-upper} since $f_3$ is constant.
\end{proof}

\begin{corollary}\label{four-dimensional-indefinite}
Let $n=3$ or $n=4$. If $S,f_3$ are constant and $S>1$, then $h$
has both a positive and a negative principal curvature at every point.
Moreover, $|X|$ is unbounded.
\end{corollary}
\begin{proof}
Suppose first that $h$ is positive definite at some point.
Then $f_3>0$.
By Proposition~\ref{mean-curvature-changes-sign}, $H<0$ at another
point. Join the two points by a path. Since principal curvatures are continuous,  
at the first point on this path where the least principal curvature
is zero, all principal curvatures are nonnegative, and at most $n-1$
are positive.
The power mean inequality gives
\[
 f_3^2=\left(\sum_i\lambda_i^3\right)^2
             \ge\frac{(\sum_i\lambda_i^2)^3}{n-1}
             =\frac{S^3}{n-1}.
\]
This contradicts Proposition~\ref{strict-cubic-obstruction}, since
$(n-2)^2\le n$ for $n=3,4$.
The same power mean inequality excludes a positive semidefinite
form with a zero eigenvalue directly.
Reverse the normal to exclude negative semidefinite forms.
Thus $h$ is indefinite everywhere.
\\
If $|X|$ were bounded, the function $|X|^2/2$ would be bounded above.
Apply Theorem~\ref{thm:Hessian} to $-|X|^2/2$.
It gives points and $\varepsilon_m\downarrow0$ at which
\[
 \nabla_i\nabla_i(|X|^2/2)\le\varepsilon_m .
\]
By \eqref{eq:psiderivs}, $\nabla_i\nabla_i(|X|^2/2)=1-H\lambda_i$.
Hence $H\lambda_i\ge1-\varepsilon_j>0$ for every $i$ and large $m$.
All principal curvatures then have the same sign, a contradiction.
Therefore $|X|$ is unbounded.
\end{proof}

\subsection{A limiting configuration at infinity}

We first record a limiting argument at points where $|X|$ tends to
infinity. The argument uses intrinsic distance.
\begin{lemma}\label{two-limits-at-infinity}
Suppose that $S>1$ is constant and $|X|$ is unbounded.
There are two limiting configurations, one with a nonpositive
limiting Hessian $M_+$ of $H$, and the other with a nonnegative
limiting Hessian $M_-$.
In each configuration a frame can be chosen so that
\begin{equation}\label{infinity-vanishing-direction}
 h_{ij}=\lambda_i\delta_{ij},\qquad
 \lambda_n=0,\qquad
 h_{ijk}=0\quad\hbox{if any index is }n,
\end{equation}
\begin{equation}\label{infinity-third-traces}
 \sum_i h_{iik}=0,\qquad \sum_i\lambda_i h_{iik}=0,
 \qquad \sum_{i,j,k}h_{ijk}^2=S(S-1).
\end{equation}
If $f_3$ is also constant, then
\begin{equation}\label{infinity-cubic-trace}
 \sum_i\lambda_i^2h_{iik}=0,
\end{equation}
and, in the respective configurations,
\begin{equation}\label{infinity-Hessian-contractions}
 \sum_i\lambda_i(M_\pm)_{ii}=S-Hf_3,
 \qquad \sum_i\lambda_i^2(M_\pm)_{ii}=f_3-Hf_4.
\end{equation}
Here all quantities in each formula refer to its own configuration.
No agreement of their spectra is asserted.
\end{lemma}
\begin{proof}
Choose $p_j$ with $|X(p_j)|\ge j^2$ since $|X|$ is unbounded.
Put $r_j(q)=d(p_j,q)$ and $\rho_j(q)=\sqrt{1+r_j(q)^2}$.
By \eqref{eq:Hbounds}, $|H|\le\sqrt{nS}$. Thus each of the functions
\[
 H-j^{-1}\rho_j,\qquad -H-j^{-1}\rho_j
\]
attains a maximum, at points denoted by $q_j^+$ and $q_j^-$.
Indeed, $\rho_j\to\infty$ at intrinsic infinity, while $H$ is bounded.
Hopf--Rinow theorem  then gives the maxima.
Comparison with the value at $p_j$ gives
\[
 r_j(q_j^\pm)\le\rho_j(q_j^\pm)
                  \le1+2j\sqrt{nS}.
\]
Consequently,
\[
 |X(q_j^\pm)|\ge j^2-1-2j\sqrt{nS}\longrightarrow\infty.
\]
The lower sectional curvature bound in \eqref{eq:curvature-lower-bounds}
and Hessian comparison give smooth upper supports for $\rho_j$
with
\[
 |\nabla\rho_j|\le1,\qquad \nabla^2\rho_j\le Cg,
\]
where $C$ is independent of $j$ and of the base point $p_j$.
Indeed, away from the cut locus,
\[
 \rho_j'=\frac{r_j}{\sqrt{1+r_j^2}},\qquad
 \rho_j''=\frac1{(1+r_j^2)^{3/2}},
\]
and these formulas, together with
$\nabla^2r_j\le(k+1/r_j)g$, give the bound.
Use the distance supports from
Lemma~\ref{trace-minimum-principle} at cut points.
At the two maxima, we therefore obtain
\[
 |\nabla H(q_j^\pm)|\le j^{-1},\qquad
 \nabla^2H(q_j^+)\le Cj^{-1}g,\qquad
 \nabla^2H(q_j^-)\ge-Cj^{-1}g.
\]
Since $|X^\perp|=|H|$ is bounded, $|X^\top(q_j^\pm)|\to\infty$.
Take $e_n=X^\top/|X^\top|$ at these points.
Equations \eqref{eq:first-position} and \eqref{mean-curvature-hessian}
give
\[
 h_{in}=\frac{H_{,i}}{|X^\top|}\longrightarrow0,
\]
and
\[
 h_{ijn}
 =\frac{H_{,ij}-h_{ij}+H\sum_kh_{ik}h_{kj}}{|X^\top|}
 \longrightarrow0.
\]
The bounds in Subsection~\ref{subsec:basic} give
\[
 |\nabla^2H|^2
 =\sum_{i,j}\left(\sum_kh_{kkij}\right)^2
 \le n\sum_{i,j,k,l}h_{ijkl}^2\le C.
\]
Thus both numerators are uniformly bounded.
\\
If necessary taking a subsequences, we can assume that  $h,\nabla h,\nabla^2h$  and
$\nabla^2H$ converge.
Rotate $e_1,\ldots,e_{n-1}$ to diagonalize the limiting $h$.
Full symmetry of $h_{ijk}$ proves
\eqref{infinity-vanishing-direction}.
Pass to the limit in
\[
 \sum_i h_{iik}=H_{,k}\longrightarrow0,\qquad
 \sum_{i,j}h_{ij}h_{ijk}=\tfrac12S_{,k}=0,\qquad
 \sum_{i,j,k}h_{ijk}^2=S(S-1).
\]
This gives \eqref{infinity-third-traces}.
If $f_3$ is constant, \eqref{eq:trace-gradients} gives
\eqref{infinity-cubic-trace}.
Finally, take limits in
\eqref{mean-curvature-moment-identities} to obtain
\eqref{infinity-Hessian-contractions}.
\end{proof}

\vskip2mm
\begin{proof}[Proof of Theorem~\ref{four-dimensional-classification}]
Suppose that $S>1$.
Theorem~\ref{thm:f3} gives $S>3$.
Corollary~\ref{four-dimensional-indefinite} shows that $|X|$ is
unbounded. Since $f_3$ is constant, we can assume  $f_3\ge0$.
Use the configuration with $M_-\ge0$ in
Lemma~\ref{two-limits-at-infinity}.

\smallskip\noindent\textit{Step 1. The limiting principal curvatures have  both signs.}
Lemma~\ref{two-limits-at-infinity} yields  $\lambda_4=0$ and
\[
 \lambda_1^2+\lambda_2^2+\lambda_3^2=S,
 \qquad
 \lambda_1^3+\lambda_2^3+\lambda_3^3=f_3.
\]
The three numbers have both signs.
Otherwise,  $f_3\ge0$  implies  that all three are nonnegative.
Thus $h\ge0$, and Cauchy--Schwarz  inequality gives
\[
 Hf_3\ge S^2,\qquad
 0\le\sum_i\lambda_i(M_-)_{ii}=S-Hf_3\le S-S^2<0,
\]
because of $S>1$.
This is impossible.
Since $h_{ij4}=0$ for any $i, j$, 
for $k=1,2,3$,  we have
\begin{equation}\label{four-dimensional-three-traces}
 \sum_{i=1}^3h_{iik}
 =\sum_{i=1}^3\lambda_i h_{iik}
 =\sum_{i=1}^3\lambda_i^2h_{iik}=0.
\end{equation}

\smallskip\noindent\textit{Step 2.  The limiting principal curvatures are distinct.}
Suppose, after permutation, that the three numbers are $a,a,b$.
They have both signs, so $ab<0$.
We can assume $a>0>b$.
This calculation uses neither the sign of $f_3$ nor the sign of $M_-$.
The first two equations in \eqref{four-dimensional-three-traces}
give
\[
 h_{33k}=0,\qquad h_{11k}+h_{22k}=0.
\]
By full symmetry of $h_{ijk}$,  we obtain
\[
 4(h_{111}^2+h_{112}^2)+6(h_{223}^2+h_{123}^2)=S(S-1),\qquad
 C=4a(h_{111}^2+h_{112}^2)+6\frac{2a+b}{3}(h_{223}^2+h_{123}^2).
\]
Since $H=2a+b$ and $a-b>0$,
\[
 C=\frac H3\{4(h_{111}^2+h_{112}^2)+6(h_{223}^2+h_{123}^2)\}+4\frac{a-b}{3}(h_{111}^2+h_{112}^2)
       \ge\frac H3S(S-1).
\]
Equation~\eqref{c3:constant-norm-C}, with $tS=S-1$, therefore implies
$3f_3\ge2SH$.
On the other hand, using $S=2a^2+b^2$ and $f_3=2a^3+b^3$, we obtain
\begin{align*}
 0\leq 3f_3-2SH
 &=-2a^3+b^3-4ab^2-4a^2b\\
 &=-a^3+b^3-a(a+2b)^2<0.
\end{align*}
This is a contradiction. Hence the three limiting principal
curvatures are distinct.
\\
The coefficient matrix in \eqref{four-dimensional-three-traces} is
\[
 \begin{pmatrix}
 1&1&1\\
 \lambda_1&\lambda_2&\lambda_3\\
 \lambda_1^2&\lambda_2^2&\lambda_3^2
 \end{pmatrix}.
\]
Its determinant is
$(\lambda_2-\lambda_1)(\lambda_3-\lambda_1)(\lambda_3-\lambda_2)\ne0$.
Therefore $h_{iik}=0$ for all $i,k\le3$.
Only the six permutations of $h_{123}$ can be nonzero.
Counting these permutations gives
\[
 C=\dfrac13\sum_{i, j, k}(\lambda_i+\lambda_j+\lambda_k)h_{ijk}^2=\frac H3\sum_{i,j,k}h_{ijk}^2
     =\frac H3S(S-1).
\]
By \eqref{c3:constant-norm-C}, we obtain
\begin{equation}\label{four-dimensional-common-H}
 H=\frac{3f_3}{2S}.
\end{equation}
If  $f_3>0$, \eqref{four-dimensional-common-H} yields  $H>0$, and
\eqref{infinity-Hessian-contractions} gives
\[
 0\le\sum_i\lambda_i^2(M_-)_{ii}=f_3-Hf_4.
\]
By \eqref{four-dimensional-common-H}, $f_4\le2S/3$.
On the other hand, Cauchy--Schwarz inequality gives
\[
 S^2=(\lambda_1^2+\lambda_2^2+\lambda_3^2)^2
 \le3(\lambda_1^4+\lambda_2^4+\lambda_3^4)=3f_4.
\]
Consequently $S\le2$, contrary to $S>3$.
\\
If  $f_3=0$, 
equation \eqref{four-dimensional-common-H} gives $H=0$.
In the configuration with $M_-\ge0$,
\eqref{infinity-Hessian-contractions} gives
\[
 \sum_i\lambda_i^2(M_-)_{ii}=0,\qquad
 \sum_i\lambda_i(M_-)_{ii}=S.
\]
Cauchy--Schwarz  inequality gives
\[
 S^2=\left(\sum_i\lambda_i(M_-)_{ii}\right)^2
 \le\left(\sum_i(M_-)_{ii}\right)
       \left(\sum_i\lambda_i^2(M_-)_{ii}\right)=0.
\]
This is impossible. 
\\
Thus, $S>1$ does not hold.  $S\leq 1$ and \eqref{simons} give $S=0$ or $S=1$.
Proposition~\ref{prop:parallel-classification} gives the stated images.
\end{proof}

\section{Proof of Theorem 1.3}
We prove Theorem~\ref{thm:maingrowth} in two steps:
\[
 \Div_\rho V=S(S-1)>0,\quad |V|\le C
 \quad\Longrightarrow\quad W_p'(r)\ge\gamma W_p(r)
 \quad\text{a.e.}
\]
Here $W_p$ is defined by \eqref{eq:ball-volumes}.

\begin{lemma}\label{prop:weightedfield}
Let $X:M^n\to\R^{n+1}$ be a self-shrinker with constant $S$ and $f_3$.
The vector field
\begin{equation}\label{eq:Vdef}
V=\sum_iV_ie_i,\qquad
V_i=\sum_j\bigl((S-1)h_{ij}-f_3\delta_{ij}\bigr)H_{,j}
\end{equation}
satisfies
\begin{equation}\label{eq:Vdiv}
\Div_\rho V=S(S-1).
\end{equation}
If $S>1$, then
\begin{equation}\label{eq:Vbound}
|V|\le \bigl((S-1)\sqrt S+|f_3|\bigr)\sqrt{nS(S-1)}.
\end{equation}
\end{lemma}
\begin{proof}
Codazzi equation and
\eqref{eq:psiderivs} give
\[
\sum_i h_{iji}=H_{,j},
\qquad \sum_i h_{ij}\psi_{,i}=H_{,j}.
\]
For any smooth function $v$, we have
\[
\begin{split}
&\Div_\rho(h\nabla v)=\Div(h\nabla v)-\langle \nabla \psi, h\nabla v\rangle\\
&=\sum_j\left(\sum_i h_{iji}
                 -\sum_i\psi_{,i} h_{ij}\right)v_{,j}
  +\sum_{i,j}h_{ij}v_{,ij}\\
&=\sum_{i,j}h_{ij}v_{,ij}.
\end{split}
\]
Take $v=H$. 
\eqref{mean-curvature-moment-identities} gives
\[
 \Div_\rho(h\nabla H)=\sum_{i,j}h_{ij}H_{,ij}=S-Hf_3.
\]
By \eqref{eq:Lh} and \eqref{eq:weighted-divergence},
$\Div_\rho(\nabla H)=(1-S)H$. Since $S$ and $f_3$ are constant,
we obtain
\begin{align*}
\Div_\rho V
&=(S-1)\Div_\rho(h\nabla H)-f_3\Lop H\\
&=(S-1)(S-Hf_3)-f_3(1-S)H\\
&=S(S-1).
\end{align*}
For any tangent vector $\vec w$, Cauchy--Schwarz inequality
\[
|h \vec w|=|\sum_i\lambda_iw_i|\le\sqrt S\,|\vec w|.
\]
Combining this inequality  with \eqref{eq:Hbounds} and
\eqref{eq:Vdef}, we obtain
\[
|V|\le\bigl((S-1)\sqrt S+|f_3|\bigr)|\nabla H|
\le\bigl((S-1)\sqrt S+|f_3|\bigr)\sqrt{nS(S-1)}.
\]
This proves \eqref{eq:Vbound}.
\end{proof}

\smallskip
We now assume $S>1$. By \eqref{eq:Vdiv} and \eqref{eq:Vbound},
we obtain the following estimate.

\begin{proposition}\label{prop:weightedgrowth}
Let $X:M^n\to\R^{n+1}$ be an $n$-dimensional complete self-shrinker with constant $S>1$ and $f_3$.
For $p\in M$, put
\[
W_p(r)=\int_{B_r^M(p)}\rho\,d\mu,
\qquad
\gamma=\frac{\sqrt{S(S-1)}}{\sqrt{n}\bigl((S-1)\sqrt S+|f_3|\bigr)}>0.
\]
Then, for $r\ge r_0>0$,
\begin{equation}\label{eq:Wgrowth}
W_p(r)\ge W_p(r_0)e^{\gamma(r-r_0)}.
\end{equation}
In particular, the Gaussian weighted volume of intrinsic balls cannot
have sub-exponential growth.
\end{proposition}
\begin{proof}
Fix $r,\delta>0$ and write $d(x)=d_M(p,x)$. Set
\[
\eta_{r,\delta}(x)=
\begin{cases}
1,&d(x)\le r,\\
(r+\delta-d(x))/\delta,&r<d(x)<r+\delta,\\
0,&d(x)\ge r+\delta.
\end{cases}
\]
By completeness, $\eta_{r,\delta}$ has compact support. Since
$|\nabla d|=1$ almost everywhere, \eqref{eq:weighted-parts}  and Lemma \ref{prop:weightedfield} give
\begin{align*}
S(S-1)\int_M\eta_{r,\delta}\rho\,d\mu
&=-\int_M\langle\nabla\eta_{r,\delta},V\rangle\rho\,d\mu\\
&=\frac1\delta\int_{\{r<d<r+\delta\}}
                  \langle\nabla d,V\rangle\rho\,d\mu\\
&\le\frac{\bigl((S-1)\sqrt S+|f_3|\bigr)\sqrt{nS(S-1)}}{\delta}
       \bigl(W_p(r+\delta)-W_p(r)\bigr).
\end{align*}
The use of the Lipschitz function $\eta_{r,\delta}$ follows by
smooth approximation on its compact support. Since
$\eta_{r,\delta}=1$ on $B_r^M(p)$, we have
\[
S(S-1)W_p(r)
\le\frac{\bigl((S-1)\sqrt S+|f_3|\bigr)\sqrt{nS(S-1)}}{\delta}
       \bigl(W_p(r+\delta)-W_p(r)\bigr).
\]
By the coarea formula, $W_p$ is locally absolutely continuous.
Letting $\delta\downarrow0$, we obtain
\[
W_p'(r)\ge\gamma W_p(r)
\qquad\text{for almost every }r>0.
\]
Since $W_p(r)>0$ for $r>0$, integration gives
\[
\log W_p(r)-\log W_p(r_0)
=\int_{r_0}^r\frac{W_p'(t)}{W_p(t)}\,dt
\ge\gamma(r-r_0).
\]
Thus \eqref{eq:Wgrowth} holds. Moreover,
\[
\liminf_{r\to\infty}\frac{\log W_p(r)}{r}\ge\gamma>0,
\]
which excludes sub-exponential growth.
\end{proof}

\smallskip
We are now in a position to prove Theorem~\ref{thm:maingrowth}.

\begin{proof}[Proof of Theorem~\ref{thm:maingrowth}]
By \eqref{simons}, $S=0$ or $S\ge1$. Proposition~\ref{prop:subexp-weighted}
shows that either growth hypothesis in the theorem implies
\eqref{eq:main-weighted-growth}. Suppose that $S>1$.
For any $r\ge r_0>0$, Proposition~\ref{prop:weightedgrowth} gives
\[
 \frac{\log W_p(r)}r
 \ge\gamma+\frac{\log W_p(r_0)-\gamma r_0}{r}.
\]
Letting $r\to\infty$, we obtain
\[
 0<\gamma\le\liminf_{r\to\infty}\frac{\log W_p(r)}r
 \le\limsup_{r\to\infty}\frac{\log W_p(r)}r=0.
\]
This is a contradiction. Hence $S=0$ or $S=1$.
The image classification follows from
Proposition~\ref{prop:parallel-classification}.
\end{proof}

\begin{remark}
The volume growth hypothesis excludes $S>1$. Without it, \eqref{eq:Wgrowth}
shows that any example with $S>1$ has volume of exponentially growing weighted
intrinsic balls.
\end{remark}

\begin{corollary}
Under constant $S>1$ and $f_3$, let $\gamma$ be as in
Proposition~\ref{prop:weightedgrowth}. Then
\[
 \liminf_{r\to\infty}\frac{\log W_p(r)}r\ge\gamma.
\]
\end{corollary}
\begin{proof}
Take logarithms in \eqref{eq:Wgrowth} and divide by $r$.
If the total Gaussian volume is finite, $W_p(r)$ is bounded above and
bounded below by $W_p(1)>0$ for $r\ge1$.
\end{proof}

 \vskip2mm
\noindent
{\bf Acknowledgement}.
This work was partly supported by MEXT Promotion of Distinctive Joint Research Center Program JPMXP0723833165 
and Osaka Metropolitan University Strategic Research Promotion Project (Development of International Research Hubs).

\end{document}